# A GENERAL THEORY OF MINIMUM ABERRATION AND ITS APPLICATIONS[1]

By Ching-Shui Cheng and Boxin Tang

*Academia Sinica and University of California, Berkeley, and Simon Fraser University and University of Memphis*

Minimum aberration is an increasingly popular criterion for comparing and assessing fractional factorial designs, and few would question its importance and usefulness nowadays. In the past decade or so, a great deal of work has been done on minimum aberration and its various extensions. This paper develops a general theory of minimum aberration based on a sound statistical principle. Our theory provides a unified framework for minimum aberration and further extends the existing work in the area. More importantly, the theory offers a systematic method that enables experimenters to derive their own aberration criteria. Our general theory also brings together two seemingly separate research areas: one on minimum aberration designs and the other on designs with requirement sets. To facilitate the design construction, we develop a complementary design theory for quite a general class of aberration criteria. As an immediate application, we present some construction results on a weak version of this class of criteria.

**1. Introduction.** The general problem considered in this paper is how to select the "best" fractional factorial designs. In situations where we have little or no knowledge about the effects that are potentially important, it is appropriate to select designs using the minimum aberration criterion [Fries and Hunter (1980)]. Wu and Hamada (2000) contains tables of many known minimum aberration designs. Minimum aberration designs enjoy some attractive robust properties [Cheng, Steinberg and Sun (1999) and Tang and Deng (1999)]. Much work has been done on the construction of minimum aberration designs. For details, we refer to Franklin (1984), Chen and Wu (1991), Chen (1992), Chen and Hedayat (1996), Tang and Wu

Received March 2003; revised June 2004.
[1]Supported by the National Science Foundation.

(1996), Suen, Chen and Wu (1997) and many others. Sitter, Chen and Feder (1997), Chen and Cheng (1999) and Cheng and Wu (2002) developed aberration criteria for blocked fractional factorials. A projective geometric approach to blocking fractional factorials is considered in Mukerjee and Wu (1999), and blocked fractional factorials with maximum estimation capacity are studied by Cheng and Mukerjee (2001). Wu and Zhu (2003) examined the use of a minimum aberration criterion for design selection in robust parameter design.

Developing a general theory of minimum aberration is motivated by the desire to unify various versions of minimum aberration that have recently appeared in the literature. Based on a sound statistical principle, this paper develops a general theory of minimum aberration and discusses its various applications. In addition to building a unified framework for many of the existing aberration criteria, the theory provides a method for deriving other aberration criteria that may be more appropriate for given design situations. A minimum aberration design can be called a model robust design because of its robust properties. A design with a requirement set [Greenfield (1976)] is a model specific design since such a design specifies a set of effects to be estimated. Our general theory is capable of bringing together these seemingly unrelated two classes of designs.

We will focus our discussion on two-level regular fractional factorial designs. However, most of our arguments are quite general. Section 2 motivates, introduces and studies a general criterion of minimum aberration and discusses its application to blocked fractional factorials, and to fractional factorials when some 2-factor interactions are important. Section 3 is devoted to developing a theory of complementary designs for quite a general class of aberration criteria, and presents some construction results on weak aberration.

In what follows, we introduce some notation and definitions to set the stage for the later development. A regular $2^{m-p}$ design has $m$ factors each at two levels and $n = 2^{m-p}$ runs, and is completely determined by $p$ independent defining words. The two levels are denoted by $+1$ and $-1$, so the design matrix $D$ of such a design is an $n \times m$ matrix of $\pm 1$. The defining relation of a $2^{m-p}$ design is the complete set of defining words. Labels of factors are referred to as letters. A defining word specifies a set of letters that has the property that the product of the corresponding columns of $D$ is a column of all plus ones. Including $I$, the column of all ones, the defining relation of a $2^{m-p}$ design has $2^p$ defining words. Let $A_i(D)$ be the number of defining words of length $i$ in the defining relation of design $D$, where the length of a word is the number of letters in the word. The resolution of design $D$ is the integer $R$ such that $A_i(D) = 0$ for $i = 1, \ldots, R-1$ and $A_R(D) > 0$. The minimum aberration criterion selects designs that sequentially minimize $A_1(D), \ldots, A_m(D)$. For designs of resolution at least III, we have $A_1 = A_2 = 0$, so the minimum aberration criterion selects designs that sequentially minimize $A_3(D), \ldots, A_m(D)$.

## 2. General theory of minimum aberration and its applications.

2.1. *A general criterion of minimum aberration.* Besides the grand mean $\gamma_0$, there are in all $2^m - 1$ factorial effects in a $2^{m-p}$ design. Suppose that out of the $2^m - 1$ effects, we are interested in estimating a set of effects $\gamma_1$. Then the fitted model is given by

$$(1) \qquad Y = \gamma_0 I + W_1 \gamma_1 + \varepsilon,$$

where $Y$ denotes the vector of $n$ observations, $\gamma_1$ the vector of the effects to be estimated, $W_1$ the model matrix corresponding to $\gamma_1$ and $\varepsilon$ the vector of uncorrelated random errors, assumed to have a zero mean and a constant variance. Because the remaining effects may not be negligible, we should choose a design that minimizes their contamination on the estimation of $\gamma_1$, from among all designs allowing estimation of the model in (1). Suppose that prior knowledge enables us to divide these remaining effects into $J - 1$ groups, denoted by $\gamma_2, \ldots, \gamma_J$, in such a way that the effects in $\gamma_j$ are more important than those in $\gamma_{j+1}$, for $j = 2, \ldots, J - 1$. Then the true model can be written as

$$(2) \qquad Y = \gamma_0 I + W_1 \gamma_1 + W_2 \gamma_2 + \cdots + W_J \gamma_J + \varepsilon,$$

where $W_j$ is the model matrix corresponding to $\gamma_j$ for $j = 1, \ldots, J$. The least-squares solution $\hat{\gamma}_1 = (W_1^T W_1)^{-1} W_1^T Y = n^{-1} W_1^T Y$ from the fitted model in (1) has expectation, taken under the true model in (2), $E(\hat{\gamma}_1) = \gamma_1 + C_2 \gamma_2 + \cdots + C_J \gamma_J$, where $C_j = n^{-1} W_1^T W_j$ for $j \geq 2$. So the bias of $\hat{\gamma}_1$ in estimating $\gamma_1$ is $C_2 \gamma_2 + \cdots + C_J \gamma_J$. Note that $C_j \gamma_j$ represents the contribution of $\gamma_j$ to the bias. As $\gamma_j$ is unknown, we will have to work with $C_j$. One size measure for a matrix $C = (c_{ij})$ is given by $\|C\|^2 \stackrel{\text{def}}{=} \text{trace}(C^T C) = \sum_{i,j} c_{ij}^2$. Since the effects in $\gamma_j$ are more important than those in $\gamma_{j+1}$, to minimize the bias of $\hat{\gamma}_1$, heuristically we can sequentially minimize $\|C_2\|^2, \ldots, \|C_J\|^2$. For regular designs, the entries of $C_j$ are either 0 or 1, and therefore $N_j = \|C_j\|^2$ is simply the number of effects in $\gamma_j$ that are aliased with those in $\gamma_1$, for $j = 2, \ldots, J$. Two effects are aliased (or confounded) with each other if their corresponding columns in the model matrix are identical.

DEFINITION 1. The general criterion of aberration is defined as the one that selects designs by sequentially minimizing $N_2, \ldots, N_J$, where $N_j$ is the number of effects in $\gamma_j$ that are aliased with those in $\gamma_1$, for $j = 2, \ldots, J$.

For convenience, the vector $(N_2, \ldots, N_J)$ is called the word length pattern with respect to $(\gamma_1, \gamma_2, \ldots, \gamma_J)$. An immediate application is to the situation where $\gamma_1$ are the main effects and $\gamma_j$ are the $j$-factor interactions. In this case, we have

$$(3) \qquad N_j = (j+1)A_{j+1} + (m-j+1)A_{j-1}$$

for $2 \leq j \leq m-1$, and $N_m = A_{m-1}$, where $A_j$ is the number of defining words of length $j$ as introduced in Section 1. The relationship in (3) leads to the conclusion that sequentially minimizing $N_2, N_3, \ldots$ is equivalent to sequentially minimizing $A_3, A_4, \ldots$.

LEMMA 1. *If $\gamma_1$ are the main effects and $\gamma_j$ are the $j$-factor interactions, then the general criterion of aberration, given in Definition 1, is equivalent to the usual criterion of aberration.*

The essential result in Lemma 1 was first given by Tang and Deng (1999), who in fact presented their result under a more general framework, where both regular and nonregular designs are considered. Superficially, Lemma 1 provides a statistical justification for the usual criterion of aberration, which was originally defined from the combinatorial point of view. A message running a bit deeper here is that the usual minimum aberration criterion of combinatorial nature can in fact be *derived* from a general theory based on a sound statistical principle.

A more general result than Lemma 1 can easily be obtained. Let $\gamma_1$ be the main effects and all the interactions involving up to $q$ factors. For the model in (1) to be estimable, a design of resolution $2q+1$ must exist, which implies that $A_i = 0$ for $i = 1, \ldots, 2q$. Now let $\gamma_j$ be the $(q-1+j)$-factor interactions for $j \geq 2$. We can easily show that

$$
\begin{aligned}
N_j = & \sum_{i=1}^{q} \binom{q-1+j+i}{i} A_{q-1+j+i} \\
(4) \qquad & + \sum_{i=1}^{q} \binom{q-3+j+i}{i-1} \binom{m-(q-3+j+i)}{1} A_{q-3+j+i} \\
& + \sum_{i=2}^{q} \binom{q-5+j+i}{i-2} \binom{m-(q-5+j+i)}{2} A_{q-5+j+i} + \cdots.
\end{aligned}
$$

Since $A_i = 0$ for $i = 1, \ldots, 2q$, we have

$$N_2 = \binom{2q+1}{q} A_{2q+1}, \qquad N_3 = \binom{2q+2}{q} A_{2q+2} + \binom{2q+1}{q-1} A_{2q+1},$$

and so on. Noting that the leading term for $N_j$ in (4) is given by $\binom{2q-1+j}{q} A_{2q-1+j}$, we conclude that sequentially minimizing $N_2, N_3, \ldots$ is equivalent to sequentially minimizing $A_{2q+1}, A_{2q+2}, \ldots$. This establishes the following result.

THEOREM 1. *If $\gamma_1$ are the main effects and all the interactions involving up to $q$ factors, and $\gamma_j$ are the $(q-1+j)$-factor interactions for $j \geq 2$, then the general criterion of aberration gives rise to the usual criterion of aberration that sequentially minimizes $A_{2q+1}, A_{2q+2}, \ldots$ among all designs of resolution $2q+1$.*

2.2. *Application to blocked fractional factorials.* In addition to $m$ treatment factors, a blocked fractional factorial contains $m_1$ blocking factors. The main effects of blocking factors are block effects. So are the interactions of blocking factors. Therefore, the total number of block effects produced by $m_1$ blocking factors is $2^{m_1} - 1$.

To avoid confusion, the terms "factor" and "effect" are carefully used in this paper. We stick to the meanings of the terms as in the following: a factor has a main effect, two factors have a 2-factor interaction (effect), three factors have a 3-factor interaction (effect) and so on. We therefore speak of $m_1$ blocking factors and $2^{m_1} - 1$ block effects.

A basic requirement for blocked fractional factorials is that all the $2^{m_1} - 1$ block effects should be included in the fitted model. In addition, interactions between treatment and blocking factors are assumed to be nonexistent, which is necessary for the effectiveness of blocking. Now consider all the treatment effects. To apply the general theory, we need to specify a set of treatment effects we want to estimate. Then the fitted model contains these treatment effects in addition to all the block effects. In what follows, we look at two important special cases.

The first case is that the main effects of the $m$ treatment factors are in the fitted model. Then $\gamma_1$ in model (1) consists of the main effects of all the $m$ treatment factors and all the $2^{m_1} - 1$ block effects. For the remaining treatment effects, we assume as usual that the hierarchical ordering principle applies [Wu and Hamada (2000)], and therefore $\gamma_j$ in model (2) represents the vector of all the $j$-factor interactions of treatment factors, where $j = 2, \ldots, m$.

A defining word in a blocked fractional factorial is a subset of $m + m_1$ letters among which $m$ letters represent treatment factors and $m_1$ letters represent blocking factors. Let $A_j$ be the number of defining words of length $j$ that contain no blocking factors, and let $B_j$ be the number of defining words that contain $j$ treatment factors and at least one blocking factor. Note that we must have $A_1 = A_2 = B_0 = B_1 = 0$ for the fitted model to be estimable.

PROPOSITION 1. *Let $\gamma_1$ denote all main effects of treatment factors and all block effects, and let $\gamma_j$ denote all the $j$-factor interactions of treatment* *factors. Then the word length pattern* $(N_2, \ldots, N_m)$ *is given by* $N_2 = 3A_3 + B_2$, $N_3 = 4A_4 + B_3$, *and in general*

$$N_j = (j+1)A_{j+1} + (m-j+1)A_{j-1} + B_j, \tag{5}$$

*where* $A_j$ *and* $B_j$ *are defined in the preceding paragraph.*

The proof is straightforward. Our general criterion of aberration for blocked fractional factorials therefore selects designs by sequentially minimizing $N_2 = 3A_3 + B_2$, $N_3 = 4A_4 + B_3$ and so on. Chen and Cheng (1999) proposed a criterion of aberration, and using our notation, their word length pattern is given by $(3A_3 + B_2, A_4, 10A_5 + B_3, A_6, \ldots)$. We see that the leading component in their criterion is identical to the leading component $N_2$ in our general criterion of aberration. Sitter, Chen and Feder (1997) also proposed a criterion that sequentially minimizes $A_3, B_2, A_4, B_3, A_5, B_4$, and so on. If the magnitude of $A_{j+1}$ is about the same as or larger than that of $B_j$, the criterion of Sitter, Chen and Feder (1997) provides a reasonably good approximation to our general criterion. We give an illustration using a simple example.

EXAMPLE 1. Suppose that we want to study nine factors in 16 runs, which are to be arranged in two blocks. We use $1, \ldots, 9$ to denote the nine factors, and $b$ to denote the single blocking factor. Consider the following two designs. Design $D_1$ is given by $5 = 123$, $6 = 124$, $7 = 134$, $8 = 234$, $9 = 12$, and $b = 13$, and design $D_2$ given by $5 = 123$, $6 = 124$, $7 = 134$, $8 = 13$, $9 = 12$, and $b = 234$. One can easily verify that $A_3(D_1) = 4$ and $B_2(D_1) = 4$, and $A_3(D_2) = 6$ and $B_2(D_2) = 2$. The criterion of Sitter, Chen and Feder (1997) selects $D_1$ as a better design because $D_1$ has a smaller value of $A_3$. Now applying our criterion, we see that $N_2(D_1) = 16$ and $N_2(D_2) = 20$, and again $D_1$ is better. Note that design $D_1$ in fact has a larger value of $B_2$ but its smaller value of $A_3$ plays a dominant role here.

The other important special case is that we are interested in estimating all main effects and all 2-factor interactions of treatment factors. So $\gamma_1$ consists of all the main effects and all the 2-factor interactions of treatment factors, as well as all the block effects. For $j \geq 2$, $\gamma_j$ is the vector of all the $(j+1)$-factor interactions of the treatment factors. For the fitted model to be estimable, we must have $A_1 = A_2 = A_3 = A_4 = B_0 = B_1 = B_2 = 0$. Applying our general theory, we obtain the following.

PROPOSITION 2. *Suppose that* $\gamma_1$ *consists of all main effects and all* 2-*factor interactions of treatment factors, as well as all block effects. Let* $\gamma_j$ *be the vector of all the* $(j+1)$-*factor interactions of the treatment factors for* $j \geq 2$. Then the word length pattern $(N_2, N_3, \ldots)$ is given by $N_2 = 10A_5 + B_3$, $N_3 = 15A_6 + 5A_5 + B_4$, and in general

$$N_j = (j+2)A_{j+2} + \binom{j+3}{2} A_{j+3} + B_{j+2}$$
(6)
$$+ (m-j)A_j + (m-j-1)(j+1)A_{j+1} + \binom{m-j+1}{2} A_{j-1},$$

where $A_j$ and $B_j$ are defined as before.

Proposition 2 is easily established by a simple combinatorial argument. Comparing our criterion with that of Chen and Cheng (1999), we find that the leading component in their criterion becomes $10A_5 + B_3$, which is precisely the $N_2$ given by our general theory. In fact, we have verified that the leading component in the word length pattern of Chen and Cheng (1999) is also correct if in addition to all block effects, the true model consists of all main effects and all interactions involving up to $q$ factors with $q \geq 3$. One can therefore appropriately regard the aberration criterion of Chen and Cheng (1999) as a robust version of our general aberration criterion when applied to blocked fractional factorials.

Before moving on, we remark that like other work in the area, block effects are treated as fixed effects in this paper. Our discussion in this section focuses on the situation where we are interested in estimating these block effects. If the block effects are not of interest, the contamination on their estimation due to nonnegligible treatment effects will not be a concern. Our general criterion can easily be modified to accommodate this situation. In the meantime, many new issues arise and they will be looked into in the future.

2.3. *Fractional factorials when some 2-factor interactions are important.* Suppose that a set of 2-factor interactions (2fi's) is postulated to be important, and in addition to the main effects, we are also interested in estimating these important 2fi's. In this situation the fitted model in (1) consists of all main effects and these important 2fi's. For the remaining effects, we assume as usual that the hierarchical ordering principle applies. Using the notation in Section 2.1, we have that $\gamma_1$ represents the main effects and the important 2fi's, $\gamma_2$ represents the remaining 2fi's and $\gamma_j$ represents the $j$-factor interactions for $j \geq 3$.

A 2fi of a fractional factorial $D$ can be represented by an unordered pair $(c, d)$, where $c$ and $d$ are two columns of $D$. Let $(c_1, d_1), \ldots, (c_S, d_S)$ denote the important 2fi's. For each 2fi $(c_s, d_s)$ where $s = 1, \ldots, S$, let $A_j(c_s, d_s)$ be the number of length-$j$ words containing both letters $c_s$ and $d_s$, let $A_j(c_s, \bar{d}_s)$ be the number of length-$j$ words containing $c_s$ but not $d_s$, let $A_j(\bar{c}_s, d_s)$ be the number of length-$j$ words containing $d_s$ but not $c_s$, and let $A_j(\bar{c}_s, \bar{d}_s)$ be the number of length-$j$ words containing neither $c_s$ nor $d_s$. Obviously, $A_j = A_j(c_s, d_s) + A_j(c_s, \bar{d}_s) + A_j(\bar{c}_s, d_s) + A_j(\bar{c}_s, \bar{d}_s)$. Let

$$A_j^{(2)} = \sum_{s=1}^{S} A_j(c_s, d_s), \qquad A_j^{(1)} = \sum_{s=1}^{S} [A_j(c_s, \bar{d}_s) + A_j(\bar{c}_s, d_s)],$$

(7)

$$A_j^{(0)} = \sum_{s=1}^{S} A_j(\bar{c}_s, \bar{d}_s).$$

If a defining word of length $j$ contains more than one pair of letters in the list of the important 2fi's $(c_1, d_1), \ldots, (c_S, d_S)$, it is counted more than once in calculating $A_j^{(2)}$. So $A_j^{(2)}$ in fact represents the total number of times that a defining word of length $j$ contains an important 2fi $(c_s, d_s)$. Interpretation of $A_j^{(1)}$ and $A_j^{(0)}$ is similar.

PROPOSITION 3. *When some 2fi's are important, the word length pattern $(N_2, N_3, \ldots, N_m)$ is given by $N_2 = 3A_3 + A_4^{(2)}$, $N_3 = 4A_4 + A_5^{(2)} + A_3^{(1)}$ and in general*

(8) $$N_j = (j+1)A_{j+1} + (m-j+1)A_{j-1} + A_{j+2}^{(2)} + A_j^{(1)} + A_{j-2}^{(0)},$$

*where $A_{j+2}^{(2)}$, $A_j^{(1)}$ and $A_{j-2}^{(0)}$ are defined in* (7).

The above version of the word length pattern is given in terms of the defining words of the original design matrix $D$. We now present another version in terms of the defining words of the augmented design given by the model matrix $W_1 = (D, D_2)$, where $W_1$ is as in model (1) and $D_2$ corresponds to the important 2fi's in the fitted model. This latter version is convenient for developing a general complementary design theory in Section 3.

Consider the words in the defining relation of the augmented design $W_1 = (D, D_2)$. Let $A_j$ be the number of length-$j$ words having all their $j$ letters from $D$, and let $B_j$ be the number of length-$(j+1)$ words having $j$ letters from $D$ and one letter from $D_2$.

PROPOSITION 4. *When some 2fi's are important, the word length pattern $(N_2, N_3, \ldots, N_m)$ is given by $N_2 = 3A_3 + B_2 - S$, and*

(9) $$N_j = (j+1)A_{j+1} + (m-j+1)A_{j-1} + B_j$$

*for $j \geq 3$.*

In Proposition 4, dependence of $N_j$ on the important 2fi's is expressed through $B_j$, which depends on matrix $D_2$, given by the columns of the important 2fi's.

The expression for $N_2$ in Proposition 4 needs a bit of explanation. Let $(c_s, d_s)$ be an important 2fi, for $s = 1, \ldots, S$. Then the three columns $c_s$, $d_s$ and $c_s d_s$, where columns $c_s$ and $d_s$ are from $D$ and column $c_s d_s$ is from $D_2$, form a word of length 3 that contributes to $B_2$ but not to $N_2$. This explains why we have $N_2 = 3A_3 + B_2 - S$ instead of $N_2 = 3A_3 + B_2$.

Ke and Tang (2003) examined practical issues in design selection using the general criterion of aberration when some 2fi's are important, and presented a collection of designs of 16 and 32 runs for models containing up to four important 2fi's.

2.4. *Other applications.* In robust parameter design, there are two sets of factors, control factors and noise factors. The goal of the experiment is to choose the settings of control factors so that the response variable is insensitive to noise factors. Suitable designs should therefore allow analysis of both location and dispersion effects. Wu and Zhu (2003) examined the use of an aberration criterion for robust parameter design which is mainly motivated by the analysis of location effects. It would be interesting to see how our general theory can be modified to take into account the analysis of dispersion effects. One possibility is to select designs using our criterion from among those designs allowing suitable analysis of dispersion effects as can be found in Hedayat and Stufken (1999). Our general theory is potentially useful in fractional factorial split plot designs. Huang, Chen and Voelkel (1998) and Bingham and Sitter (1999) considered aberration criteria for split plot designs. Since split plot designs have more than one error structure, some sort of modification seems necessary for our theory to be applicable to such problems. In our future research, both areas of application will be considered.

**3. Theory of complementary designs.** In this section, we will develop a complementary design theory for a class of aberration criteria. This class of criteria, to be introduced below, is quite broad, and in particular includes as special cases the aberration criteria for blocked fractional factorials and for designs when some 2fi's are important as discussed in Sections 2.2 and 2.3, respectively.

3.1. *A class of aberration criteria.* Suppose that besides the main effects $\gamma_1'$ of $m$ factors, we are also interested in estimating additional $S$ effects $\gamma_1''$. For convenience, these $m$ factors are called major factors. In addition to the $m$ major factors, we may have $m_1$ minor factors, where $m_1 \geq 0$. When $m_1 = 0$, the additional $S$ effects $\gamma_1''$ are a set of interactions only involving major factors. When $m_1 \geq 1$, the effects in $\gamma_1''$ are a subset of effects from the collection of all the following effects: the interactions only involving major factors, the main effects of minor factors, the interactions only involving minor factors and the interactions involving both major and minor factors. The $\gamma_1$ in the fitted model (1) is therefore given by $\gamma_1 = (\gamma_1', \gamma_1'')$. Let $\gamma_j$ denote the $j$-factor interactions only involving the major factors that are not included in $\gamma_1''$, for $j = 2, \ldots, m$. We assume as earlier that the effects in $\gamma_j$ are more important than those in $\gamma_{j+1}$, for $j \geq 2$. With $\gamma_1, \ldots, \gamma_m$ defined above, the true model is now given in (2). A remark on the true model is in order when there is at least one minor factor, that is, $m_1 \geq 1$. An implicit assumption made here is that all other effects involving at least one minor factor besides those in $\gamma_1''$ are assumed to be nonexistent. The above formulation is fairly general, and includes as special cases all the situations discussed in Sections 2.1–2.3. For example, for blocked fractional factorials, we take the treatment factors as the major factors and the blocking factors as the minor factors. Choices for major and minor factors are also natural for fractional factorials in the row-column setting [Cheng and Mukerjee (2003)].

We now derive the word length pattern for the above situation. Let $D_1$ be the design matrix corresponding to the main effects $\gamma_1'$ of the $m$ major factors and let $D_2$ be the matrix corresponding to the additional $S$ effects $\gamma_1''$. Note that $D_1$ has $m$ columns and $D_2$ has $S$ columns. The word length pattern will be given in terms of the model matrix $W = (D_1, D_2)$, specified by its two components $D_1$ and $D_2$. Let $A_j(D_1)$ be the number of length-$j$ defining words in design $D_1$. Define $B_j(D_1, D_2)$ to be the number of length-$(j+1)$ defining words in design $W = (D_1, D_2)$, which have $j$ letters from $D_1$ and one letter from $D_2$. Then it is easily established that the word length pattern $(N_2, \ldots, N_m)$ is given by

$$N_j(D_1, D_2) = (j+1)A_{j+1}(D_1) + (m-j+1)A_{j-1}(D_1) + B_j(D_1, D_2) - S_j$$

for $j \geq 2$, where $S_j$ is the number of the interactions of $j$ major factors that are included in $\gamma_1''$. (For an explanation of why $S_j$ is necessary, see the end of Section 2.3.) Note that $S_j$ is a constant for the purpose of choosing $D_1$ and $D_2$. For simplicity, ignoring $S_j$, we redefine the word length pattern $(N_2, \ldots, N_m)$ as

(10) $\quad N_j(D_1, D_2) = (j+1)A_{j+1}(D_1) + (m-j+1)A_{j-1}(D_1) + B_j(D_1, D_2).$

The goal here is to choose $D_1$ and $D_2$ by sequentially minimizing $N_2, N_3, \ldots$.

3.2. *A complementary design theory.* Let $H_k = (D_1, D_2, D_3)$, where $H_k$ denotes a saturated design of $n = 2^k$ runs and $n - 1$ factors. Obviously, it is impossible to completely characterize design pair $(D_1, D_2)$ through $D_3$ alone. Our complementary design theory to be developed below characterizes design pair $(D_1, D_2)$ through design pair $(D_2, D_3)$. This approach is most effective when the number of columns in $D_3$ is smaller than that in $D_1$.

We need a result from Tang and Wu (1996) and Suen, Chen and Wu (1997), who developed a complementary design theory for the usual minimum aberration criterion. The explicit coefficients in Lemma 2 are due to Suen, Chen and Wu (1997).

LEMMA 2. *Let $H_k = (D, \bar{D})$, where $D$ has $m$ factors. Then we have*

$$A_j(D) = \sum_{i=0}^{j} c_m(i,j) A_i(\bar{D}),$$

*where $c_m(1,j) = c_m(2,j) = 0$, $c_m(i,j) = (-1)^{j-[(j-i)/2]} \binom{m-2^{k-1}}{[(j-i)/2]}$ for $3 \leq i \leq j$, and*

$$c_m(0,j) = (-1)^{j-[j/2]} \binom{m-2^{k-1}}{[j/2]} + 2^{-k}[P_j(0;m) - P_j(2^{k-1};m)],$$

*where $P_j(x;m) = \sum_{s=0}^{j}(-1)^s \binom{x}{s}\binom{m-x}{j-s}$ is a Krawtchouk polynomial.*

Note that $N_j(D_1, D_2)$ in (10) depends on design pair $(D_1, D_2)$. The main result of our complementary design theory is contained in the following theorem, which expresses $N_j(D_1, D_2)$ in terms of design pair $(D_2, D_3)$.

THEOREM 2. *The word length pattern in (10) for the class of criteria discussed in Section 3.1 depends on design pair $(D_2, D_3)$ through*

$$N_j(D_1, D_2) = \sum_{i=0}^{j+1}[(j+1-S)c_m(i,j+1) + Sc_{m+1}(i,j+1)]A_i(D_2 \cup D_3)$$

$$- \sum_{i=0}^{j+1} c_{m+1}(i,j+1) E_i(D_2, D_3)$$

$$+ (m-j+1)\sum_{i=0}^{j-1} c_m(i,j-1) A_i(D_2 \cup D_3),$$

*where $E_i(D_2, D_3) = \sum_{p=1}^{i} p E_i^{(p)}(D_2, D_3)$ with $E_i^{(p)}(D_2, D_3)$ denoting the number of length-$i$ defining words in $D_2 \cup D_3$ that have exactly $p$ letters from $D_2$.*

PROOF. Applying Lemma 2 to design $D_1$, we have

$$(11) \qquad A_j(D_1) = \sum_{i=0}^{j} c_m(i,j) A_i(D_2 \cup D_3).$$

For any $d \in D_2$, applying Lemma 2 to design $D_1 \cup \{d\}$, we obtain

$$(12) \qquad A_j(D_1 \cup \{d\}) = \sum_{i=0}^{j} c_{m+1}(i,j) A_i((D_2 \setminus \{d\}) \cup D_3).$$

Let $B_{j-1}(d, D_1)$ be the number of length-$j$ defining words in design $(D_1, D_2)$ that contain letter $d$ and $j-1$ letters from $D_1$. Clearly, we have $A_j(D_1 \cup \{d\}) = A_j(D_1) + B_{j-1}(d, D_1)$. Let

$$(13) \qquad T_j(d, D_2, D_3) = \sum_{i=0}^{j} c_{m+1}(i, j) A_i((D_2 \setminus \{d\}) \cup D_3).$$

Then (12) can be rewritten as

$$(14) \qquad A_j(D_1) + B_{j-1}(d, D_1) = T_j(d, D_2, D_3).$$

Taking summation over all $d$ in $D_2$ on both sides of (14) gives

$$\sum_{d \in D_2} (A_j(D_1) + B_{j-1}(d, D_1)) = \sum_{d \in D_2} T_j(d, D_2, D_3).$$

Noting that $D_2$ has $S$ columns and that $B_{j-1}(D_1, D_2)$ defined in Section 3.1 is equal to $\sum_{d \in D_2} B_{j-1}(d, D_1)$, we obtain $SA_j(D_1) + B_{j-1}(D_1, D_2) = T_j(D_2, D_3)$, where

$$(15) \qquad T_j(D_2, D_3) = \sum_{d \in D_2} T_j(d, D_2, D_3).$$

Therefore, $B_j(D_1, D_2) = T_{j+1}(D_2, D_3) - SA_{j+1}(D_1)$. Substituting this expression of $B_j(D_1, D_2)$ into (10), we obtain

$$(16) \quad \begin{aligned} &N_j(D_1, D_2) \\ &= (j+1-S)A_{j+1}(D_1) + T_{j+1}(D_2, D_3) + (m-j+1)A_{j-1}(D_1). \end{aligned}$$

Now let us calculate $T_{j+1}(D_2, D_3)$ in (15). Let $E_i(d, D_2, D_3)$ be the number of length-$i$ defining words in $(D_2, D_3)$ that contain letter $d$. We have $A_i((D_2 \setminus \{d\}) \cup D_3) = A_i(D_2 \cup D_3) - E_i(d, D_2, D_3)$. Then (13) becomes

$$T_{j+1}(d, D_2, D_3) = \sum_{i=0}^{j+1} c_{m+1}(i, j+1)[A_i(D_2 \cup D_3) - E_i(d, D_2, D_3)].$$

Summing both sides over all $d$ in $D_2$, we obtain

$$(17) \quad \begin{aligned} T_{j+1}(D_2, D_3) &= S \sum_{i=0}^{j+1} c_{m+1}(i, j+1) A_i(D_2 \cup D_3) \\ &\quad - \sum_{i=0}^{j+1} c_{m+1}(i, j+1) E_i(D_2, D_3), \end{aligned}$$

where $E_i(D_2, D_3) = \sum_{d \in D_2} E_i(d, D_2, D_3)$. Note that $E_i(d, D_2, D_3)$ is the number of length-$i$ words containing letter $d$ in design $D_2 \cup D_3$. Thus $E_i(D_2, D_3)$ represents the total number of times length-$i$ words in design $D_2 \cup D_3$ contain a letter in $D_2$. Therefore

$$(18) \qquad E_i(D_2, D_3) = \sum_{p=1}^{i} p E_i^{(p)}(D_2, D_3),$$

where $E_i^{(p)}$ denotes the number of length-$i$ words in $D_2 \cup D_3$ having exactly $p$ letters from $D_2$. Combining (11) and (16)–(18), we obtain the result in Theorem 2. □

Chen and Cheng (1999) developed a complementary design theory for blocked fractional factorials. Our complementary design theory given in Theorem 2 is applicable to a broad class of aberration criteria including all the cases discussed in Sections 2.1–2.3. We want to mention that our theory does not include their theory as a special case, because our word length pattern when applied to blocked factorials is not exactly the same as theirs. On the other hand, one can adopt our approach to derive their complementary design theory. Our approach appears considerably simpler than theirs.

Zhu (2003) found a relationship between $A_{ij0}$ and $A_{0kl}$, where $A_{ij0}$ is the number of length-$(i+j)$ words having $i$ letters from $D_1$ and $j$ letters from $D_2$, and $A_{0kl}$ is the number of length-$(k+l)$ words having $k$ letters from $D_2$ and $l$ letters from $D_3$. In principle, Theorem 2 is derivable from his result. On the other hand, it is not obvious how one can obtain the result in Theorem 2, which clearly shows how $N_2$ depends on $D_2$ and $D_3$, from Zhu's rather involved formula that connects $A_{ij0}$ to $A_{0kl}$.

3.3. *Some results on weak aberration.* The general criterion of aberration sequentially minimizes $N_2(D_1, D_2), N_3(D_1, D_2), \ldots$. A weak version of the criterion is given by minimizing $N_2(D_1, D_2) = 3A_3(D_1) + B_2(D_1, D_2)$ alone. Using Theorem 2, we find that

$$N_2(D_1, D_2) = \text{constant} - 3A_3(D_2 \cup D_3) + E_3(D_2, D_3),$$

where the *constant* does not depend on $D_2$ and $D_3$. Noting that

$$E_3(D_2, D_3) = E_3^{(1)}(D_2, D_3) + 2E_3^{(2)}(D_2, D_3) + 3E_3^{(3)}(D_2, D_3),$$

$$A_3(D_2 \cup D_3) = A_3(D_3) + E_3^{(1)}(D_2, D_3) + E_3^{(2)}(D_2, D_3) + E_3^{(3)}(D_2, D_3),$$

we have that

$$N_2(D_1, D_2) = \text{constant} - 3A_3(D_3) - 2E_3^{(1)}(D_2, D_3) - E_3^{(2)}(D_2, D_3).$$

So minimizing $N_2(D_1, D_2)$ is equivalent to maximizing

$$(19) \qquad g(D_2, D_3) = 3A_3(D_3) + 2E_3^{(1)}(D_2, D_3) + E_3^{(2)}(D_2, D_3).$$

The following lemma gives an upper bound on $g(D_2, D_3)$.

LEMMA 3. *Let $m_2 = S$ and $m_3 = 2^k - 1 - m - S$ be the numbers of columns in $D_2$ and $D_3$, respectively. We have that:*

(i)  $3A_3(D_3) + E_3^{(1)}(D_2, D_3) \leq \binom{m_3}{2}$,
(ii) $E_3^{(1)}(D_2, D_3) + E_3^{(2)}(D_2, D_3) \leq m_2 m_3 / 2$, *and*
(iii) $g(D_2, D_3) \leq m_3(m_2 + m_3 - 1)/2$.

*The upper bound in* (iii) *is reached if and only if the bounds in* (i) *and* (ii) *are both reached.*

PROOF. For any two columns $c$ and $d$ in $D_3$, the product $cd$ must belong to one of $D_1$, $D_2$ or $D_3$. Consider all the $\binom{m_3}{2}$ pairs of columns in $D_3$. The number of the pairs whose products are in $D_2 \cup D_3$ is given by $E_3^{(1)}(D_2, D_3) + 3A_3(D_3)$. Therefore

$$E_3^{(1)}(D_2, D_3) + 3A_3(D_3) \leq \binom{m_3}{2}, \tag{20}$$

which proves part (i) of Lemma 3. Similarly, by considering all the products $cd$ such that $cd \in D_2 \cup D_3$ where $c \in D_2$ and $d \in D_3$, we obtain

$$2E_3^{(1)}(D_2, D_3) + 2E_3^{(2)}(D_2, D_3) \leq m_2 m_3, \tag{21}$$

from which Lemma 3(ii) follows. Combining (20) and (21), we obtain

$$g(D_2, D_3) \leq \binom{m_3}{2} + m_2 m_3 / 2 = m_3(m_2 + m_3 - 1)/2.$$

This is Lemma 3(iii). The last statement in Lemma 3 is obvious. □

From the proof of Lemma 3, we see that the bound in (20) is reached if $cd \in D_2 \cup D_3$ for any two columns $c, d \in D_3$, and that the bound in (21) is reached if $cd \in D_2 \cup D_3$ for any $c \in D_2$ and any $d \in D_3$. One structure for $D_2$ and $D_3$ to have these properties is given as follows. Let $a_1, a_2, \ldots, a_k$ be a set of $k$ independent columns that generates the saturated design $H_k$ of $n = 2^k$ runs and $n - 1$ factors. Now choose $D_2 \cup D_3$ to be $H_r = H(a_1, \ldots, a_r)$, the saturated design generated by independent columns $a_1, \ldots, a_r$ where $r = 1, \ldots, k - 1$. Note that $D_2$ can be any $m_2 = S$ columns from $H_r$.

THEOREM 3. *Let $H_r$ be the saturated design generated by $r$ independent columns $a_1, \ldots, a_r$. Then so long as $D_2 \cup D_3 = H_r$, any design pair $(D_2, D_3)$ maximizes $g(D_2, D_3)$ in* (19). *Therefore design pair $(D_1, D_2)$ has minimum weak aberration, where $D_1$ is given by $H_k \setminus H_r$.*

Recall that in introducing the class of aberration criteria in Section 3.1, design $D_1$ corresponds to the main effects of $m$ major factors and design $D_2$ represents $m_2 = S$ additional effects we are interested in estimating, which may involve some minor factors. In order for design pair $(D_1, D_2)$ given in Theorem 3 to be a legitimate design, we need to specify $D_2$ in such a way that it indeed represents the $S$ additional effects. We now look at two situations. The first is that $D_1$ represents the main effects of $m$ treatment factors, and $D_2$ are all the $2^{m_1} - 1$ block effects given by $m_1$ blocking factors. Then choosing $D_2 = H_{m_1}$, the saturated design generated by $a_1, \ldots, a_{m_1}$ where $m_1 \leq r$ in Theorem 3, satisfies the requirement. This characterization for blocked designs was given in Chen and Cheng (1999). We see that it is now derived from Theorem 3. The second situation we will look at is that $D_1$ are the main effects of $m$ factors, and $D_2$ are some 2-factor interactions of the $m$ factors. We illustrate how to choose $D_2$ through an example.

EXAMPLE 2. Suppose that we want a 16-run design that allows estimation of the main effects $1, 2, 3, 4, 5, 6, 7$ and 8 of eight factors and the following 2-factor interactions: $12, 13, 24$ and $35$. Let $a_1, a_2, a_3, a_4$ be four independent columns. Theorem 3 says that we should choose $D_2 \cup D_3 = \{a_1, a_2, a_1a_2, a_3, a_3a_1, a_3a_2, a_3a_1a_2\}$. That is, the eight factors are assigned to the columns in $D_1 = \{a_4, a_4a_1, a_4a_2, a_4a_1a_2, a_4a_3, a_4a_3a_1, a_4a_3a_2, a_4a_3a_1a_2\}$. Now assign factor 1 to $a_4$, factor 2 to $a_4a_1$, factor 3 to $a_4a_2$, factor 4 to $a_4a_3$ and factor 5 to $a_4a_3a_1$. Factors $6, 7$ and 8 can be arbitrarily assigned to the remaining three columns in $D_1$. We have that $12 = a_1$, $13 = a_2$, $24 = a_1a_3$ and $35 = a_1a_2a_3$. So $D_2 = \{a_1, a_2, a_1a_3, a_1a_2a_3\}$.

**Acknowledgments.** The authors thank an Associate Editor and two referees for constructive comments.

Department of Statistics
University of California, Berkeley
Berkeley, California 94720-3860
USA
e-mail: cheng@stat.berkeley.edu

Department of Statistics
and Actuarial Science
Simon Fraser University
Burnaby, British Columbia
Canada V5A 1S6
e-mail: boxint@cs.sfu.ca